\documentclass[12pt]{article}

\usepackage{amsmath, amsthm, amssymb}

\newtheorem*{thm*}{Theorem}
\newtheorem{thm}{Theorem}

\newtheorem{lem}{Lemma}

\newtheorem{prop}{Proposition}
\newcommand{\C}{\mathbb C}
\newcommand{\R}{\mathbb R}

\newcommand{\N}{\mathbb N}

\newcommand{\f}{\displaystyle \frac}
\newcommand{\Z}{\mathbb Z}
\theoremstyle{definition}

\theoremstyle{remark}
\newtheorem{rem}{Remark}

\DeclareSymbolFont{AMSb}{U}{msb}{m}{n}
\DeclareMathSymbol{\bbbn}{\mathalpha}{AMSb}{"4E}
\DeclareMathSymbol{\bbbz}{\mathalpha}{AMSb}{"5A}
\DeclareMathSymbol{\bbbr}{\mathalpha}{AMSb}{"52}
\DeclareMathSymbol{\bbbq}{\mathalpha}{AMSb}{"51}
\DeclareMathSymbol{\bbbc}{\mathalpha}{AMSb}{"43}

\newcommand{\imag}{\mathrm{i}}
\newcommand{\Id}{\mathrm{Id}}
\newcommand{\const}{\mathrm{const}}
\newcommand{\Chi}{\Psi}
\newcommand{\geoint}{G}
\newcommand{\diag}{\mathop{\mathrm{diag}}}
\newcommand{\trace}{\mathop{\mathrm{tr}}}
\newcommand{\Spectr}{\mathop{\mathrm{Spectr}}}

\numberwithin{equation}{section}

\title{On the analytic non-integrability   of the Rattleback problem }

\author{H. R. Dullin$^*$, A. V. Tsygvintsev$^\dagger$ \\ \\\parbox{8cm}{\begin{center}
$^*$ Department of Mathematical Sciences\\ Loughborough University\\ Loughborough\\
LE11 3TU, UK \\{\tt H.R.Dullin@lboro.ac.uk}
\end{center}}\hfill\parbox{8cm}{\begin{center}
$^\dagger$ Unit\'e de math\'ematiques pures et appliqu\'ees\\  Ecole Normale Sup\'erieure de Lyon\\ 46, all\'ee d'Italie, Lyon\\
F--69364  Lyon Cedex 07, France\\{\tt atsygvin@umpa.ens-lyon.fr}\end{center}}}

\begin{document}

\bibliographystyle{amsplain}

\maketitle

\begin{abstract}
\noindent We establish the  analytic non-integrability of the nonholonomic ellipsoidal 
rattleback model for a large class of parameter values.
Our approach is based on the study of the monodromy group of the normal variational equations around a particular orbit. The embedding of  the equations of the heavy rigid body  into the rattleback model is discussed. 
\end{abstract}

\section{Introduction}

The rattleback's amazing mechanical behaviour can be described as follows: when spun on a flat horizontal surface in the clockwise  direction this top  continues to spin in the same direction until it consumes the initial spin energy;  when, however, spun in the counterclockwise direction, the spinning  soon ceases, the body briefly oscillates, and then reverses its spin direction and thus spins in the clockwise direction until the energy is consumed.  The first mathematical model of this phenomena belongs to Walker (1896) who studied the linearized equations of motion and  concluded  that the completely stable motion is possible in only one  (say clockwise)  spin direction.  This  classical explanation of the rattleback's behavior is incomplete since it does not reflect the global dynamical effects explaining   the transfer of trajectories from the vicinity of the unstable solution to the stable one.  To analyze thoroughly this question we propose to study the  nonholonomic equations of the rattleback in the complex domain   and particularly to study \emph{the existence of  additional analytic  first integrals}. It has been observed in  many mechanical systems that non existence of analytic  first integrals  is usually associated to complicated chaotic behavior of trajectories of the system.
 We mention that actually only numerical evidence for chaos  in the rattleback systems has been observed \cite{Bor}.

In our case the rattleback represents a  full ellipsoidal  body whose  center of mass $P$   coincides with  its  geometric one. All vectors are defined in the body fixed axes (with origin in $P$) coinciding the with principal geometric axes of the ellipsoid
\begin{equation}\label{el}
E(r) = 
\frac{r_1^2}{b_1^2}+\frac{r_2^2}{b_2^2}+\frac{r_3^2}{b_3^2}
=1 \,.
\end{equation}

Let $\langle \,,\, \rangle $ denote the euclidean scalar product in $\R^3$ and let $||\cdot||$ be  the corresponding  norm.  In addition, let $[\,, \,]$ denote the cross product in $\R^3$.

As shown  in \cite{Garcia},  the nonholonomic equations of motion  can be written in the following form
\begin{subequations}
\begin{align} \label{1eqs}
\Theta\, \f{d\omega}{dt}+m\,[r,[\f{d\omega}{dt},r]] &=
-[\omega,\Theta \omega]
-m\, [r,[\omega,[\omega,r]]]  
+ mg\, [r,\gamma]-m\, [r,[\omega,\f{dr}{dt}]]\,,
\\ \label{2eqs}
\f{d\gamma}{dt} &=[\gamma,\omega]\,,
\end{align}
\end{subequations}
where 
\begin{equation}\label{eqs}
\Theta=\left[ \begin {array}{ccc} \Sigma_{{11}}&\Sigma_{{12}}&0
\\\noalign{\medskip}\Sigma_{{12}}&\Sigma_{{22}}&0
\\\noalign{\medskip}0&0&\Sigma_{{33}}\end {array} \right]\,,
\end{equation}
is the inertia matrix; $m$ -- mass of the body; $g$ -- gravitational constant;  $\omega=(\omega_1,\omega_2,\omega_3)^T$ -- angular velocity; $\gamma=(\gamma_1,\gamma_2,\gamma_3)^T$  -- unit vector normal to the body's surface at the point of contact; $r=(r_1,r_2,r_3)^T$-- position of contact point.

To observe the rattleback behavior,  one  has to  impose the following condition on $\Theta$ 
(see  e.g.~\cite{Bor}) and $b_i$
\begin{equation}
\begin{array}{llll}
\left[ \begin {array}{cc} \Sigma_{{11}}&\Sigma_{{12}}
\\\noalign{\medskip}\Sigma_{{12}}&\Sigma_{{22}}\end {array} \right] 
= \left[ \begin {array}{cc} I_{{1}}  \cos^2 \delta  
 +I_{{2}} \sin^2 \delta  &
 \left( I_{{1}}-I_{{2}} \right) \cos \delta  \sin
 \delta  \\\noalign{\medskip} \left( I_{{1}}-I_{{2}}
 \right) \cos \delta  \sin \delta  &I_{{1}
} \sin^2 \delta   +I_{{2}} \cos^2
 \delta   \end {array} \right], \\  \\ \Sigma_{33}=I_3, \quad b_1\neq b_2\,,
\end{array}
\end{equation}
where $I_1$, $I_2$, $I_3$ are the principal components of the inertia tensor  whose  principal horizontal axes  are rotated by the angle $\delta >0$ with respect to the $r_1$-$r_2$-axes of the ellipsoid.

Unless otherwise mentioned below we will  assume that
\begin{equation} \label{condition2}
\quad \Sigma_{ij}>0, \quad \Sigma_{3,3}<\Sigma_{1,1}+\Sigma_{2,2}, \quad b_1> b_2>0, \quad b_3>0,\quad m>0, \quad g>0\,.
\end{equation}

Using \eqref{el},  the vector  $r=R(\gamma)$ is  defined by solving $\gamma = - \nabla E/||\nabla E||$
for $r$, which gives

\begin{equation}\label{rg}
 r_{{i}}=R_i(\gamma) = \displaystyle\frac {-b_i^2\gamma_i}
 {\sqrt {b_1^2\gamma_1^2+b_2^2\gamma_2^2+b_3^2\gamma_3^2}}, \quad i = 1,2,3 \,.
\end{equation} 

Solving   \eqref{1eqs} for $\dot \omega$ and using \eqref{rg} we transform the nonholonomic rattleback equations  to the standard form
\begin{equation} \label{EQS}
\begin{array}{lll}
\f{ d\omega}{dt}=F(\omega,\gamma),\\ \\
\f{ d\gamma}{dt}=[\gamma,\omega]\,, 
\end{array}
\end{equation} 
where $F=(F_1,F_2,F_3)^T$ is a vector field rational in the variables $\omega$, $\gamma$ and $s={\sqrt {b_1^2\gamma_1^2+b_2^2\gamma_2^2+b_3^2\gamma_3^2}}$. 

 Let $(\omega,\gamma,s)\in \C^7$. Then the vector field \eqref{EQS}, as a function of $(\omega,\gamma,s)$, is  analytic in the domain $\mathcal{D}=\C^7\setminus (\mathcal{S}\cup \mathcal{L})$ where $\mathcal{L}=\{(\omega,\gamma,s)\in \C^7\,:\, s=0\}$ and  $\mathcal{S}\subset \C^7$ is the surface on which the determinant of the matrix $U$ given below is zero
\begin{equation}
U=
\Theta - m [r, [r, \cdot]] = \Theta + m \big( \langle  r,r\rangle  ) \Id - r \otimes r\big) \,.
\end{equation}
This matrix appears when we solve  \eqref{1eqs} to find  $\dot \omega$.
For an arbitrary non-zero vector $u \in \R^3$ 
\begin{align*}
\langle  u, U u\rangle &
= \langle  u, \Theta u\rangle   + m \langle  u,u\rangle   \langle  r,r\rangle   - m \langle  u,r\rangle  ^2 \\
& = \langle  u, \Theta u\rangle   +m || u \times r||^2 >  0 \,,
\end{align*} 
so that $U$ is  positive  definite and hence $\det(U)=0$ never occurs in the mechanical case.

The  system \eqref{EQS} always possesses two first integrals  (see \cite{Garcia}):
\begin{equation} \label{i1}
H=\f{m}{2}||[\omega,R(\gamma)]||^2+\f{1}{2}\langle  \omega,\Theta \omega\rangle  -mg\langle  R(\gamma),\gamma\rangle  =h, \quad h \in \R \,,
\end{equation} 
-- the energy;
\begin{equation} \label{i2}
\geoint =\langle  \gamma,\gamma\rangle  =\gamma_1^2+\gamma_2^2+\gamma_3^2=l, \quad l \in \R_+\,,
\end{equation}
-- the geometric integral.

We see that  $H(\omega,\gamma,s)$ (after introduction of  $s$  with the help of  \eqref{rg}) and $\geoint (\gamma)$ are  analytic  functions of
 $\omega$, $\gamma$, $s$ in $\mathcal{D}$. The natural question arises  whether the  equations  \eqref{EQS} can  have a third first integral analytic (meromorphic)  in  $(\omega,\gamma,s)\in \C^7$ and functionally independent of $H$ and $\geoint $. 
We remark that the absence of first integrals in the meromorphic case is generally harder to 
prove than in the analytic one.

The paper is organized as follows.
In Chapter \ref{invariant_manifold} we describe one particular solution of the Rattleback problem. Chapter \ref{monodromy_group} discusses some properties of the normal variational equations and the monodromy group generators of this solution.  
Chapter \ref{log_branching_case} examines the case in which  the normal variations equations have logarithmic branching points. Under these conditions we show the non-existence of additional meromorphic first integrals (Theorem  \ref{main:lem}). Finally, Chapter   
\ref{analytic_first_integrals} contains the proof of  our main theorem about the analytic non-integrability:

\begin{thm} \label{MainTh} Assume that  the conditions \eqref{condition2} and   \eqref{CONDITION} hold. Then  the energy $H$ and  the geometric integral $\geoint $ are the only  
integrals of the rattleback problem \eqref{EQS} which are
analytic in $(\omega,\gamma,s)$.
\end{thm}

\subsection*{Rigid body limiting case}

The rattleback equations \eqref{1eqs}-\eqref{2eqs} formally contain the equations of the heavy rigid body in the
singular limiting case $m\to 0$, $mg \not \to 0$. Moreover the functions $R_i(\gamma)$ are
constants denoted by $r_i$ again, now designating the position of the center of mass in the body frame.
By a rotation the tensor of inertia $\Theta$ can be diagonalised so that  $\Sigma_{12}$ can 
be set to zero in this case, thus violating \eqref{condition2}. 
Under these assumptions the system is Hamiltonian with another well known integral
\begin{equation} \label{4int}
    L = \langle  \gamma, \Theta \omega\rangle   = \gamma_1 \omega_1 I_1 + \gamma_2 \omega_2 I_2 + \gamma_3 \omega_3 I_3 \,.
\end{equation}

\section{The invariant manifold \label{invariant_manifold}}

 We consider the vector field \eqref{EQS} as a function of six variables $(\omega,\gamma)$. Let $l>0$ (in the mechanical case we can set $l=1$).
For $(\omega,\gamma)\in \R^6$  the square root  in $s$ is assumed  to be always positive.
It is straightforward  to check that  the system \eqref{EQS} has the invariant manifold $M=\{(\omega,\gamma)\in \R^6\, : \, \omega_1=\omega_2=0,\gamma_3=0\}$. The mechanical sense of the motion   on $M$ is quite clear: it corresponds to rolling (or oscillating) of the body on the line of intersection of its surface with $r_3=0$ plane. We note that this  invariant manifold  exists  because of the  assumption that the $r_3$  ellipsoidal axis coincides with one of the principal inertia body axes.

Complexifying  \eqref{EQS}  we denote by $U_{l,h} \subset M\subset  \C^6$ its  orbit  in $M$ 
corresponding to fixed $(l,h)\in \C^2$ defined as the intersection of the two algebraic surfaces 
\begin{subequations}
\begin{align}
{\omega_{{3}}}^{2} \left( {\frac {m \left(b_1^4\gamma_1^2+b_2^4\gamma_2^2 \right) }{b_1^2\gamma_1^2+b_2^2\gamma_2^2}}+\Sigma_{{33}} \right) +2\,
{\it mg}\,\sqrt {b_1^2\gamma_1^2+b_2^2\gamma_2^2}& =2\,h \in \R,\\
\gamma_1^2+\gamma_2^2& =l \in \R\,.
\end{align}
\end{subequations}
We assume now  that $h\neq0$ and $l=0$. Obviously the solution obtained in this way cannot be real;
it is non-mechanical.

Let $\Gamma\subset \C^6=(\omega,\gamma)$ be the complex curve defined by
\begin{equation}\label{G}
 \omega_3=p,\quad \gamma_1=\f{p^2-\alpha^2}{\beta},\, \gamma_2=i\f{p^2-\alpha^2}{\beta},\,\omega_1=\omega_2=\gamma_3=0, \,  p\in \C \setminus \{-\alpha,+\alpha \}\,,
\end{equation}
 where  
\begin{equation}
\alpha^2=\f{2h}{ m(b_1^2 + b_2^2)+\Sigma_{3,3} }, \quad
 \beta=-\f{2mg\sqrt{b_1^2-b_2^2}}{ m(b_1^2+b_2^2)+\Sigma_{3,3}}, \quad h\in \C\,,
\end{equation}
This is the parametrization of $U_{0,h}$.

Let $P_0=(\omega_0,\gamma_0)\in \Gamma$  such that $\gamma_0\neq 0$.

The function $s=\sqrt{ b_1^2\gamma_1^2+b_2^2\gamma_2^2+b_3^2\gamma_3^2}$, once the 
square-root branching is fixed,  is analytic in a small neighborhood $U_{P_0}\subset \C^6$ of   $P_0$. We can analytically continue $s$ and all its derivatives along  $\Gamma$ with the help of the parametrization above: $s((\omega,\gamma)\in\Gamma)=s(p)=\sigma \f{p^2-\alpha^2}{\beta}$ where $\sigma=\sqrt{b_1^2-b_2^2}>0$. In particular, it shows that $s$ is  single-valued on $\Gamma$ and hence  is analytic in a small neighborhood $B_{\Gamma}\subset \C^6$ of it. Thus, the vector field \eqref{EQS}  restricted to $B_{\Gamma}$  is   an analytic function of the complex variables $(\omega,\gamma)$ and has  the invariant curve $\Gamma$.

We use the reparametrisation 
$b_i^2 = \rho_i(\rho_1-\rho_2)$, with $\rho_1 > \rho_2$,  $\rho_i = \const > 0$,
which makes $\sigma = \rho_1 - \rho_2$ polynomial on the solutions we are considering.
In particular for this solution $r = (-\rho_1, -\imag \rho_2, 0)^T$, which is a constant.

Since in the {\em rigid body limiting case} $r$ is not proportional to $\gamma$ the invariant 
manifold $M$ only exists under the additional assumption that $r_3 = 0$. 

\section{The variational equations and its monodromy group \label{monodromy_group}}

The relation between the parameters $p$ and  $t$ can be easily found by substitution of the parametrization \eqref{G} in one of the equations \eqref{EQS}. This gives 
\begin{equation} \label{changet}
dp=  i\f{p^2-\alpha^2}{2}dt\,.
\end{equation}

\begin{rem}
The change of time \eqref{changet} has infinitely many sheets. Nevertheless, one can always replace
the vector field $F(X)$, $X=(\omega,\gamma)^T$  defined in \eqref{EQS} by $ \tilde  F=(F/(2^{-1} i (p^2-\alpha^2)),1)^T$ i.e.\ consider the new autonomous  system of differential equations
\begin{equation} \label{tEQS}
\f{dX}{d\tau}=\f{\tilde F(X)}{2^{-1} i (p^2-\alpha^2)},\quad \f{dp}{d\tau}=1\,,
\end{equation}
which  will have  the same particular solution $\Gamma$ defined by \eqref{G} where $p$ is replaced by $\tau$. 

Obviously,  an autonomous analytic  first integral of  \eqref{EQS}  gives a first integral of the same type for  the vector field \eqref{tEQS}.  One sees also that $\tilde F$  is analytic in the neighborhood of $\Gamma$. The further analysis based on the variational equations of \eqref{EQS} or \eqref{tEQS} will be essentially the same.
\end{rem}

Let $\tilde \omega_i$, $\tilde \gamma_i$, $i=1,2,3$ be the variations of variables $\omega_i$, 
$\gamma_i$, respectively.
We  use the following notation for  coordinates of the variation vector  $V=(\tilde X,\tilde Y)^T=((\tilde \omega_1,\tilde \omega_2,\tilde \gamma_3),(\tilde \gamma_1,\tilde \gamma_2,\tilde \omega_3))^T$.
In these variables, the variational equations of \eqref{EQS} along $\Gamma$, after the substitution 
\eqref{changet}, have the block diagonal form 
\begin{equation} \label{inNVE}
\f{dX}{dp}= \left[ \begin {array}{cc} M&O\\\noalign{\medskip}O&N\end {array}
 \right] X, \quad X\in T_{\Gamma}B_{\Gamma}\,.
\end{equation}
where $M(p)$, $N(p)\in GL(3,\C(p))$ and $O$ is  the zero    $3\times 3$  matrix.

One can show that the derivatives of the first integrals $H$ and $\geoint $  with respect to the variables $\omega_1$, $\omega_2$, $\gamma_3$ vanishes along $\Gamma$.  So, the  linear first integrals $\tilde H =\langle  dH(\Gamma),\tilde Y\rangle  $ and $\tilde \geoint =\langle  d\geoint ,\tilde Y\rangle  $  are not useful in solving the first block system
\begin{equation} \label{NVE}
\f{d\tilde X}{dp}=M \tilde X\,,
\end{equation}
whereas  the second subsystem $\f{d\tilde Y}{dp}=N \tilde Y$  can be completely solved in radicals with the help of $\tilde H$ and $\tilde \geoint $.

As seen from the parametrization \eqref{G}, the equations \eqref{NVE} are the normal variational equations (see 
\cite{Zig} for definition) of the system \eqref{EQS} along the orbit $\Gamma$. The following lemma holds.

\begin{lem}[{\cite{Morales}},{\cite{Zig}}]  \label{MainProp}
Let us suppose that the system \eqref{EQS} has a third first integral $H_3(\omega,\gamma)$ which is  analytic (meromorphic)   in the neighborhood $B_{\Gamma}$ and functionally independent of $H$ and $\geoint $. Then the monodromy and the differential Galois groups of  the normal variational equations $\eqref{NVE}$  have a non trivial polynomial (rational) invariant.
\end{lem}

We can put the linear system \eqref{NVE} into the Fuchsian form with help of the rational transformation
\begin{equation}
\tilde X=\diag(p,p,p^2-\alpha^2) x\,,
\end{equation}
where $x=(x_1,x_2,x_3)^T\in \C^3$. The new system takes the following form
\begin{equation} \label{NVE1}
\f{dx}{dp}=T(p)x=\left(\f{A}{p-\alpha}+\f{A}{p+\alpha}+\f{B}{p}\right)x, \quad p \in \mathcal{R}=\bar \C \setminus \{-\alpha,+\alpha,0,\infty\}\,.
\end{equation}
Clearly, the monodromy groups of the systems \eqref{NVE} and \eqref{NVE1} are equivalent.
The constant $3\times3$ matrices $A$, $B$  are given in Appendix A.

\begin{rem}
One sees that the above linear equations are invariant under the transformation $p\mapsto -p$.
That in particular allows to reduce the number of finite singularities up to two  via  introducing the new time $p=\tau^2$. Nevertheless,  we prefer to keep this symmetry inside the equations  and to use it later in the study of the monodromy group of \eqref{NVE1}. 

\end{rem}

\begin{prop} \label{propos3}
The general solution of \eqref{NVE1}  is meromorphic in a neighborhood of $p=0$. 
\end{prop}
 
 The proof  follows from the fact that the matrix $M(p)$ in  \eqref{NVE} is holomorphic in the neighborhood of $p=0$  as seen  from \eqref{G}, \eqref{1eqs}. In particualr,  we have 
 \begin{equation} \label{spectrB}
 \Spectr(B)=\{0,-1,-1\}\,.
\end{equation}

That can be verified directly with help of formulas for $B$ given in Appendix A.

We will need some technical results.  The following proposition will play a crucial role for  our non-integrability  results below.
\begin{prop}\label{propos1}
The characteristic polynomials $P(x)$, $P_{\infty}(x)$  of the residue matrices  $A$  at $p=\pm \alpha$ and $A_{\infty}=-2A-B$  at $p=\infty$  always have non-real roots. 
\end{prop}
\begin{proof}

A direct calculation shows that  $P(x)$ is of the form 
\begin{equation}\label{CHARAC}
P(x)=x^3+x^2+\theta_1 x+\theta_0\,, 
\end{equation}
where with the definitions
\begin{equation}
\Chi_n = \begin{pmatrix}
\Sigma_{11} & \Sigma_{12} \\ 
\Sigma_{12} & \Sigma_{22}
\end{pmatrix}, 
\quad
\Chi_d = \Chi_n - (\Sigma_{11} + \Sigma_{22} - \Sigma_{33}) \Id, \quad
v =  \begin{pmatrix} r_1 \\ r_2 \end{pmatrix} =  \begin{pmatrix} -\rho_1 \\  -\imag \rho_2 \end{pmatrix} \,,
\end{equation}
the coefficients of $P$ are given by 
\begin{subequations} \label{coefficientsP}
\begin{align}
   \theta_d &= \det \Chi_d + m \langle  v, \Chi_d v\rangle   \\
   \theta_d \theta_1 &= \det \Chi_n + m\{ \langle   v, \Chi_n v\rangle   + \rho_3( \Sigma_{11} \rho_1 - \Sigma_{22} \rho_2 + \imag \Sigma_{12} (\rho_1 + \rho_2) ) \} \\
  \theta_d \theta_0 &= \theta_d \theta_1 - m \rho_3(\rho_1 - \rho_2) (\Sigma_{11} + \Sigma_{22} - \Sigma_{33})  \,.
\end{align}
\end{subequations}
Therefore
\begin{equation}
\mathrm{Im} \, \f{1}{\theta_0-\theta_1}=
\frac{2 \Sigma_{12} \rho_1 \rho_2}{\rho_3(\rho_1 - \rho_2) (\Sigma_{11} + \Sigma_{22} - \Sigma_{33})}\,,
\end{equation}
which is obviously non-zero according to \eqref{condition2}.
A direct calculation shows that $P_\infty(x)$ is of the form
\begin{equation}
P_\infty(x) = x^3 - 4 x^2 + \chi_1 x + \chi_0 
\end{equation}
In the notation of the previous proposition and the Appendix the linear coefficient is given by 
\begin{equation}
 \chi_1 = 5 + 4 a_1^2 + 4 a_2 a_3 + 4(c_3 - \imag c_1)/\beta
\end{equation}
The coefficients of $P(x)$ and $P_\infty(x)$ are related by 
\begin{equation} \label{coefrelation}
12 (\theta_0 - \theta_1) + \chi_0 + 3 \chi_1 = -9 \,.
\end{equation}
Therefore
\begin{equation}
\mathrm{Im} \, \f{1}{\chi_0 + 3 \chi_1 + 9}=
-\frac{ \Sigma_{12} \rho_1 \rho_2}{6 \rho_3(\rho_1 - \rho_2) (\Sigma_{11} + \Sigma_{22} - \Sigma_{33})}
\end{equation}
\end{proof}

In the rigid body case the matrices $A$ and $B$ become
\[
A = \begin{pmatrix}
0 & \imag(\Sigma_{33}-\Sigma_{22})/\Sigma_{11} &  0\
\imag(\Sigma_{11}-\Sigma_{33})/\Sigma_{22} & 0 & 0 \\
-1/\beta & -\imag/\beta & -1 
\end{pmatrix}
\]
\[
B = \begin{pmatrix}
-1 & 0 & -2 \imag mg r_2 / \Sigma_{11} \\
0 & -1 & 2 \imag mg r_1/ \Sigma_{22} \\
0 & 0 & 0 
\end{pmatrix}\,.
\]
Therefore {\small
\begin{equation*} \label{eh}
\Spectr(A) = \left \{ -1, 
- \sqrt{ \frac{(\Sigma_{11} - \Sigma_{33})(\Sigma_{22} - \Sigma_{33}) }{\Sigma_{11} \Sigma_{22} } }, 
 \sqrt{ \frac{(\Sigma_{11} - \Sigma_{33})(\Sigma_{22} - \Sigma_{33}) }{\Sigma_{11} \Sigma_{22} } } 
\right \}.
\end{equation*} }
For the Kovalesvkaya case $\Sigma_{11} = \Sigma_{22}$, $\Sigma_{33} = 2 \Sigma_{11}$
the eigenvalues become $\{ -1, 1, 1\}$ and $\{-2, 3, 3\}$ at infinity with a non-trivial Jordan block.

The following proposition shows that the   system  \eqref{NVE1} is not solvable in the Lappo-Danilevsky sense (see \cite{Er}).

\begin{prop} \label{propos2}
The residue matrices $A$ and $B $ do not have common eigenvectors. In particular, they  do not commute. 
\end{prop}

\begin{proof}
The proof is a straightforward, but since the proposition will play  an important role later, we indicate some basic steps of it.
A direct computation gives for the eigenvectors $\mathcal{V}_i$ and the corresponding eigenvalues 
$\lambda_i$  of $B$: 
$$
\begin{array}{lll}
\mathcal{V}_1=(1,v_1,v_2)^T, & \lambda_1=0,\\
\mathcal{V}_2=(1,0,0)^T, & \lambda_1=-1,\\
\mathcal{V}_3=(0,1,0)^T, & \lambda_1=-1\,,
\end{array}
$$ 
where $v_i$ are some known expressions. 

One verifies that the matrix $R=(v_1,v_2,v_3)$ is not singular once the condition \eqref{condition2} holds.  Under the same conditions  it is easy to show that $\mathcal{V}_2$ and $\mathcal{V}_3$  are not   eigenvectors of $A$. Here is one method  to prove the same statement for $\mathcal{V}_1$. We consider the conjugation 
$\tilde A = R^{-1} A R$ of $A$ by  $R$.  It is sufficient to show then that $\tilde A_{1,2}\neq 0$ or  $\tilde A_{1,3}\neq 0$.   That can be done quite easily once the condition \eqref{condition2} is  fixed. 
Otherwise,  if $C=[A,B]$,  then the non-commutativity  property follows  directly from   $C_{3,1}=1/ \beta\neq 0$.

\end{proof}

Fixing a basepoint $e\in \mathcal{R}$,  one defines the monodromy group $G$ of the system  \eqref{NVE1} as the image of   the fundamental group $\pi(\mathcal{R},e)$
given by the analytical continuations of the fundamental matrix solution $\Sigma(p)$, $\Sigma(e)=\Id$ of \eqref{NVE1} along all closed paths $\Gamma\in \pi(\mathcal{R},e)$ starting from $e$.  One verifies that    $\mathrm{tr}(A),\mathrm{tr}(B)\in \mathbb{Z}$. Hence,  $G \subset SL(3,\C)$.  According to Proposition  \ref{propos3}, the group $G=\langle M_+,M_- \rangle$ is generated by  the local monodromy transformations $M_+$, $M_-$ around   the singularities  $p=\alpha$ and $p=-\alpha$, respectively.

\section{On the non-existence of  additional meromorphic first integrals in the logarithmic branching case \label{log_branching_case}}

The next proposition will show that the logarithmic branching of solutions of the normal variational equations \eqref{NVE1} is not compatible with existence  of a new \emph{meromorphic} first integral of the rattleback problem \eqref{EQS}. 
\begin{thm} \label{main:lem}
Let us assume that the general solution of the variational equations \eqref{NVE1} has  logarithmic branching around one of the singularities $p=\pm \alpha$.  Then  the rattleback equations do not have any new meromorphic first integrals.
\end{thm}

\begin{proof}
By a suitable linear transormation   of  the general solution of \eqref{NVE1}, we can always reduce  one of the monodromy matrices around $p=\alpha$ or $p=-\alpha$  to its Jordan form. In the case that the general solution has logarithmic branching points, at least one of the monodromy matrices $M_+$, $M_-$ has a non-trivial Jordan block.

Since $p=\pm \alpha$ enter in \eqref{NVE1} in a symmetric way, it is sufficient to consider the case of  $p=\alpha$. Firstly, we  assume that $M_+$ is of the  form  

\begin{equation} \label{fm}
M_+=\left[ \begin {array}{ccc} q&1&0\\\noalign{\medskip}0&q&1
\\\noalign{\medskip}0&0&q\end {array} \right] \,.
\end{equation}
Since $M_+\in SL(3,\C)$,  one has  $q^3=1$. 

Let 
\begin{equation}
\Spectr(A)=(\lambda_1,\lambda_2,\lambda_3)\,.
\end{equation}
As follows  from \eqref{NVE1} (see e.g.~\cite{G}) 
 \begin{equation} \label{sprel}
 \Spectr(M_+)=\Spectr(M_-)=\{e^{2\pi i \lambda_1},e^{2\pi i \lambda_2},e^{2\pi i \lambda_3}\}\,.
 \end{equation} 
 
  Thus,  in  the case \eqref{fm} all eigenvalues of $A$ must be real (even rational) numbers. This is impossible according to Proposition \ref{propos1}. \\

We suppose  now that  $M_+$ is of the form

\begin{equation} \label{mon1}
M_+=\left[ \begin {array}{ccc} u&1&0\\\noalign{\medskip}0&u&0
\\\noalign{\medskip}0&0&u^{-2}\end {array} \right], \quad u\in \C^*\,.
\end{equation}

In this case we prove the following lemma.

\begin{lem} \label{LEMMA1}
Let $R(x)$, $x=(x_1,x_2,x_3)^T$ be a rational invariant of $G$.  Then $R=R(x_2,x_3)$.
\end{lem}

\begin{proof}
We  write $R(x)$, $\deg R=m$ as follows
\begin{equation}
R=\f{ \sum\limits _{i=0}^l\, x_3^{l-i}\, P_i(x_1,x_2) }  { \sum\limits _{j=0}^k \, x_3^{k-j}  Q_j(x_1,x_2) },\quad l,k \in \{0,1,2, \dots \}, \quad m=l-k\,,
\end{equation}
where $P_i(x_1,x_2)$, $Q_j(x_1,x_2) $  are homogeneous polynomials of degrees $i$ and $j$, respectively.
Let us assume that at least one of the polynomials $P_i$ or  $Q_j$  depends on $x_1$. So, for example, we will find    $0\leq \rho, s \leq l$ such that  the term $x_3^{\rho} P_s(x_1,x_2)$, $\partial P/ \partial x_1\neq 0$, $\mathrm{deg} (P_s)=s$ is a  semi-invariant of $M_+$, i.e.~it satisfies 
\begin{equation} \label{proport}
\tilde M_+^n \, P_s(x_1,x_2)=c^n\, P_s(x_1,x_2), \quad c \in \C^*,\quad \forall\, n \in \mathbb N, \quad 
\tilde M_+=\left[ \begin {array}{cc} u&1\\\noalign{\medskip}0&u
\end {array}
 \right] \,.
\end{equation} 
Since $P_s(x_1,x_2)$ is homogeneous of degree $s$, we can put it in the form 
\begin{equation}
P_s(x_1,x_2)=\alpha_0 x_1^s+\alpha_1x_1^{s-1}x_2+\cdots+\alpha_s x_2^s\,.
\end{equation}

We observe that $\tilde M_+^n (x_1,x_2)^T=(u^nx_1+nu^{n-1}x_2,u^nx_2)^T$. 
Consequently, the polynomial $\tilde M_+^n \,P_s(x_1,x_2)$ will have its coefficient of $x_2^s$ equal to 
$\theta_n=\sum _{r=0}^s\, \alpha_rn^{s-r}u^{s(n-1)+r}$.  Let $q$ be the smallest index such that $\alpha_q \neq 0$. Then the  coefficient before  $x_1^q x_2^{s-q}$  in  $\tilde M_+^n \,P_s(x_1,x_2)$  is  $\eta_n=u^{nq}\alpha_q$. Clearly, as $n\to \infty$, the asymptotic behaviors  of $\theta_n$ and $\eta_n$ are different.  So, the equality \eqref{proport} cannot be true. 
That finishes the proof of Lemma \ref{LEMMA1}.
\end{proof}

\begin{lem} \label{LL2}
All rational homogeneous invariants of $M_+$ are functions of  the invariant  $x_2^2x_3$.
\end{lem}

The proof follows easily from the condition that $\Spectr(A)$ is non-real.

\begin{lem} \label{form_int}
Let us assume that the monodromy group $G=\langle M_-,M_+\rangle$   has a rational homogeneous  invariant $R(x)$. Then  $x_2^2x_3$ is the invariant of $G$ and  $R=(x_2^2x_3)^m$ for a certain $m\in \Z^*$.
\end{lem}

\begin{proof}
Let $R(x_1,x_2,x_3)$ be a rational homogeneous invariant of $G$ written as
\begin{equation}
R=\f{ \sum\limits_{i_1+i_2+i_3=M} a_{i_1,i_2,i_3}\, x_1^{i_1}x_2^{i_2}x_3^{i_3} }{ \sum\limits_{j_1+j_2+j_3=N} b_{j_1,j_2,j_3}\, x_1^{j_1}x_2^{j_2}x_3^{j_3}    }, \quad i_k, j_p\in \{0,1,2,\dots \}\,,
\end{equation}
where $M,N$ are non-negative integer numbers and all coefficients $a_i$, $b_i$ are different from zero.

Firstly, we assume that $R$ reduces to the single term

\begin{equation} \label{1term}
R=x_1^{l_1}x_2^{l_2}x_3^{l_3}, \quad l_i\in \Z, \quad(l_1,l_2,l_3)\neq (0,0,0)\,.
\end{equation} 

Then $l_1=0$ according to Lemma \ref{LEMMA1}, and, as follows from   Lemma \ref{LL2}, there exists $m\in \Z^*$ such that 
\begin{equation} \label{Rform}
R=(x_2^2x_3)^m\,,
\end{equation}
that in turn proves  the result.

If $R$ is not of  the form \eqref{1term}, it can be written as follows
\begin{equation}\label{eqf}
R=\f{ax_1^{r_1}x_2^{r_2}x_3^{r_3}+bx_1^{q_1}x_2^{q_2}x_3^{q_3}+\cdots}{cx_1^{p_1}x_2^{p_2}x_3^{p_3}+\cdots  }, \quad r_i,q_j,p_k\in\{0,1,2,\dots\}\,,
\end{equation}
where $a,b,c\neq 0$ and $(r_1,r_2,r_3)\neq(q_1,q_2,q_3)$, $\sum r_i=\sum q_i=M$.

We note that, as seen from \eqref{eqf},  the  terms  $x_1^{r_1}x_2^{r_2}x_3^{r_3}$ and $x_1^{q_1}x_2^{q_2}x_3^{q_3}$ are   multiplied by the same constant under the action of $M_+$. Therefore, the division by $x_1^{q_1}x_2^{q_2}x_3^{q_3}$ shows that $r(x)=x_1^{k_1}x_2^{k_2}x_3^{k_3}$, $k_i=r_i-q_i$, $(k_1,k_2,k_3)\neq (0,0,0)$, $\sum k_i=0$  is the invariant of $M_+$. We have $k_1=0$ according to Lemma \ref{LL2} and hence $r(x)=(x_2 / x_3)^{k_2}$, $k_2\neq 0$. This  is  impossible according to  Lemma \ref{LL2}.

Thus, we can assume $R$ to have the form \eqref{Rform}.
This immediately fixes the  monodromy transformation $M_-$, also preserving $R$,  as follows
\begin{equation} 
M_-=\left[ \begin {array}{ccc} u_{{1}}&n&k
\\\noalign{\medskip}0&u_{{2}}&0\\\noalign{\medskip}0&0&u_{{3}}
\end {array} \right]
, \quad u_1u_2u_3=1, \quad u_i, n,k \in \C\,.
\end{equation}
We will show that  $u_1=u_2=u$ and $u_3=u^{-2}$.

 Indeed,  since $\Spectr(M_+)=\Spectr(M_-)=\{u,u,u^{-2}\}$,   we have necessarily $u_1=u$ or $u_2=u$. Let  $u_1=u$, $u_2=u^{-2}$, $u_3=u$. With help of  $R$ we obtain $u^{-3m}=1$ and so $u$ is the root of unity. Let now $u_1=u^{-2}$, $u_2=u$, $u_3=u$. Then $u^{3m}=1$ and we conclude as above.
Hence, the monodromy transformation $M_-$ takes  the form

\begin{equation} \label{mon2}
M_-=\left[ \begin {array}{ccc} u&n&k
\\\noalign{\medskip}0&u&0\\\noalign{\medskip}0&0&u^{-2}
\end {array} \right]\,.
\end{equation}
The rational function $R=x_2^2x_3$ is clearly  invariant under the action of  $M_-$ and $M_+$. 
\end{proof}

We shall follow the approach close  to that of Tannakian  (see e.g.~\cite{Morales}).

Let $\Sigma(p)$ be the fundamental matrix solution of  the  equations \eqref{NVE1} and $\Sigma^{-1}(p)=(\Sigma_1,\Sigma_2,\Sigma_3)^T$ where $\Sigma_i$  are linearly independent vector functions.  It is known (see e.g.~\cite{T1}, p. 246) that if $R(x)=x_2^2x_3$ is a polynomial invariant of the monodromy group $G$ 
then 
\begin{equation} \label{nofactor}
I(p,x)=R(\Sigma(p)^{-1}x)=\langle  \Sigma_2(p),x\rangle  ^2\langle  \Sigma_3(p),x\rangle  \,,
\end{equation}
 will be a first integral of \eqref{NVE1} invariant under the action of $G$, i.e.\ single-valued as a function of $p\in \mathcal R$. One can express this fact by stating that all  coefficients $a_{i,j,k}(p)$ in the expression for $I$,
\begin{equation} \label{u}
I=\sum_{i+j+k=3} a_{i,j,k}(p) \,x_1^i x_2^j x_3^k\,,
\end{equation}
 are rational functions of  $p$.

Since the polynomial  $I$ given by \eqref{u} can  be factorized as  \eqref{nofactor},  there exist two $3$-vector functions  $\mathcal A(p)$, $\mathcal B(p)$, algebraically dependent on $p$, such that $I$ becomes
\begin{equation} \label{fint}
I(p,x)=\langle  \mathcal A(p),x\rangle  ^2\langle  \mathcal B(p),x\rangle   \,.
\end{equation}

The  structure of the monodromy group generated by the transformations \eqref{mon1}, \eqref{mon2} suggests that the system \eqref{NVE1}  has two linearly independent particular solutions $X_{1,2}(p)$ of the following  form

\begin{equation} \label{sol1}
X_1(p)=(p-\alpha)^{\lambda_1}(p+\alpha)^{ \lambda_1}\, R_1(p),  \quad e^{2\pi i \lambda_1}=u, \quad \lambda_1\in \Spectr(A)\,,
\end{equation}

\begin{equation} \label{sol2}
X_2(p)=(p-\alpha)^{\lambda_1}(p+\alpha)^{ \lambda_1}R_2(p)+u^{-1}X_1(p)\left ( \f{\log(p-\alpha)}{2\pi i} +n\f{\log(p+\alpha)}{  2 \pi i  }     \right)\,,
\end{equation}
where $R_{1,2}(p)$ are  vector functions  rationally dependent on $p$.  
\begin{lem} \label{zeros}
We have $\langle  \mathcal A,R_1\rangle  =\langle  \mathcal B,R_1\rangle  =0$.
\end{lem} 
\begin{proof}
Plugging   the solution \eqref{sol1} into the integral \eqref{fint} we obtain 
$$ I(p, X_1)=(p-\alpha_1)^
{3\lambda_1}(p+\alpha)^{3 \lambda_1}\langle  \mathcal A,R_1\rangle  ^2\langle  \mathcal B,R_1\rangle  =c=\const\,.
$$ 
In the case $c\neq 0$, since $\mathcal A,\mathcal B,R_1$ are algebraic with respect to $p$, the last equality  implies   $ \lambda_1\in \mathbb Q$.   According to  \eqref{mon1}
\begin{equation}\label{sdvig}
\Spectr(A)=\{\lambda,\lambda+k, -2\lambda-k-1  \}\,, 
\end{equation}
for   certain $\lambda\in \C$ and $k\in \Z$.  

So, it is clear  that if $\lambda_1\in \R$, then $\Spectr(A)$ has to be  real.  According to Proposition \ref{propos1}  we therefore get $\langle  \mathcal A,R_1\rangle  ^2\langle  \mathcal B,R_1\rangle  =0$ for all $p$.  

We shall  consider the case  $\langle  \mathcal A,R_1\rangle  =0$ and $\langle  \mathcal B,R_1\rangle  \neq 0$. Obviously, the shift  $\tilde I(p,x)=I(p,x+X_1(p))$ of the first integral $I$   is again a first integral of the system \eqref{NVE1} and 
\begin{equation}\label{tttt}
\tilde I = I+(p-\alpha)^{\lambda_1}(p+\alpha)^{ \lambda_1}\langle  \mathcal B,R_1\rangle  \langle  \mathcal A,x\rangle  ^2=\const\,.
\end{equation}
Since the  equations \eqref{NVE1} are homogeneous, with help of \eqref{tttt}, \eqref{fint}, we derive  the    first integral 

\begin{equation} \label{J1}
J_1=\f{\langle  \mathcal B,x\rangle  }{(p-\alpha)^{\lambda_1}(p+\alpha)^{\lambda_1}(\mathcal B,R_1)}={(p-\alpha)^{-\lambda_1}(p+\alpha)^{-\lambda_1}}\langle  \tilde B,x\rangle  \,,
\end{equation}
where $\tilde B$ is algebraic on $p$.

Finally, combining $J_1$ with $I$,  one  gets 

\begin{equation} \label{J2}
\begin{array}{lll}
J_2= \sqrt{ \langle  \mathcal B,R_1\rangle  (p-\alpha)^{\lambda_1}(p+\alpha)^{\lambda_1}}\langle  \mathcal A,x\rangle  =\\ \\={(p-\alpha)^{\lambda_1/2}(p+\alpha)^{\lambda_1/2}}\langle  \tilde A,x\rangle  \,,
\end{array}
\end{equation}

-- the  first integral of \eqref{NVE1} with  $\tilde A(p)$ algebraic on $p$.

Obviously, $J_1$ and $J_2$ are functionally independent.  Indeed,  the vectors $\mathcal A(p)$, $\mathcal B(p)$ are independent as being  proportional to the lines $\Sigma_{2,3}(p)$ of the matrix $\Sigma^{-1}(p)$ whose determinant is not   identically zero.
Finally, in order to find the third linear first integral, we  apply the Liouville theorem to the fundamental matrix solution of \eqref{NVE1} formed by the columns $X_1(p)$, $X_2(p)$ and  the arbitrary  solution $x(p)$:
\begin{equation} \label{Sig}
\begin{array}{lll}
\Sigma(p)=(X_1(p),X_2(p),x), \\ \\  \det(\Sigma)=(p-\alpha)^{2\lambda_1}(p+\alpha)^{2\lambda_1}
\det(R_1,R_2,x)=\const \cdot \, a(p)\,,
\end{array}
\end{equation}
where $a(p)$ is a rational  function of $p$ in view of  $\trace(A),\trace(B)\in \Z$.

One derives from \eqref{Sig} the  following first integral of \eqref{NVE1}
\begin{equation}\label{J3}
J_3=(p-\alpha)^{-2\lambda_1}(p+\alpha)^{-2\lambda_1}\langle  \mathcal C,x\rangle  \,,
\end{equation}
with $\mathcal C(p)$ algebraic on $p$.

In the case the vectors $\tilde A,\tilde B,\mathcal C$ are linearly independent,  finding $x$ from the linear system $J_1=c_1,J_2=c_2,J_3=c_3$, $c_{1,2,3}\in \C$, we see that the general solution of \eqref{NVE1} does  not contain  logarithmic branching   and so we get  a contradiction. Let us assume  that $\tilde A,\tilde B,\mathcal C$ are linearly dependent.  Then   $\mathcal C=l_1\tilde A+l_2\tilde B$ where   $l_1$, $l_2$ are certain  algebraic functions of $p$ and the following relation  holds  in view of \eqref{J1}, \eqref{J2}, \eqref{J3}
\begin{equation} \label{relation_c_i}
(p-\alpha)^{-\frac{5}{2}\lambda_1}(p+\alpha)^{-\frac{5}{2}\lambda_1}c_1l_1+(p-\alpha)^{-\lambda_1}(p+\alpha)^{-\lambda_1}c_2l_2=c_3\,.
\end{equation}

Since $J_1$, $J_2$ are functionally independent, the last expression  shows that $\lambda_1\in \mathbb Q$   and the proof  is finished with help of \eqref{sdvig}. The  case $\langle  \mathcal A,R_1\rangle  \neq 0$, $\langle  \mathcal B,R_1\rangle  =0$ is treated in the analogous way.
\end{proof}

With help of the Lemma~\ref{zeros} one obtains  $I(p,x+X_2)=I(p, x+(p-\alpha)^{\lambda_1}(p+\alpha)^{\lambda_1}R_2)$ which is a first integral of \eqref{NVE1}.   As before we  show  $\langle  \mathcal A,R_2\rangle  =\langle  \mathcal B,R_2\rangle  =0$. Since $\mathcal A$, $\mathcal B$  are linearly independent, it follows from \eqref{sol1}, \eqref{sol2} that  $X_2=\theta X_1$ for a certain function $\theta=\theta(p)\neq \const$.  This situation can be ruled out by the substitution of $X_1$, $X_2$  into the system 
\eqref{NVE1} that gives a  contradiction.  According to Lemma \ref{MainProp} the proof of Proposition 
\ref{main:lem} is finished.

\end{proof}

\section{Non-existence of additional analytic first integrals  \label{analytic_first_integrals}}

Our aim now is to prove the non-existence of new analytic first integrals of the rattleback problem \eqref{EQS} under one of the following  hyperbolicity  conditions  
\begin{equation} \label{condition-1}
\Spectr(M_{\pm})=\{s_{1,2,3}\in  \C\, :\, 0<|s_1|<1,  |s_2|>1, |s_3|>1, s_2\neq s_3 \}\,,
\end{equation}
or
 \begin{equation} \label{condition-2}
\Spectr(M_{\pm})=\{s_{1,2,3}\in  \C\, :\, |s_1|>1, 0< |s_2|<1, 0<|s_3|<1, s_2 \neq s_3  \}\,,
 \end{equation}
 which represent  certain restrictions on the  eigenvalues of the characteristic polynomial \eqref{CHARAC}.
Recalling  that   $\sum_i \lambda_i=-1$ and $s_1 s_2 s_3=1$ we deduce from  the relations  $e^{2\pi i \lambda_i}=s_i$ that  at least one of the  conditions  \eqref{condition-1} or  \eqref{condition-2} is satisfied when 
\begin{equation} \label{CONDITION}
\lambda_1, \lambda_2, \lambda_3\not \in \R \quad \mathrm{and} \quad \mathrm{Im}\, \lambda_i -\mathrm{Im}\, \lambda_j \neq 0, \quad \forall\, i\neq j\,,
\end{equation}
where $\lambda_i$, $i=1,2,3$ are three  roots of the cubic algebraic equation
\begin{equation}
x^3+x^2+\theta_1 x +\theta_0=0\,,
\end{equation}
whose coefficients depend  on ($\Sigma_{i,j}$, $b_i$, $m$) and  are defined  by  \eqref{coefficientsP}.

From now on, we will assume that the condition \eqref{CONDITION}  is fulfilled  and  that the case \eqref{condition-1} holds.

 Let $\mathcal L$ be the space of linear forms $l=\langle  L,x\rangle  $, $L\in \C^n$ dual to $\C^n$. To each $M\in GL(n,\C)$ we associate the linear automorphism $M\, :\, \mathcal L\to \mathcal L$ according to $M\,.\,l=\langle  M^TL,x\rangle  $. 
The next results shows how the hyperbolicity of $G$  implies the reducibility of its polynomial invariants. 

\begin{prop} \label{pi}
Let the monodromy group $G=\langle M_+,M_-\rangle$ have a polynomial homogeneous invariant $P(x_1,x,_2,x_3)$. Then  one of the following situations holds:

\noindent a) $P=x_1^{\rho}x_2^l L(x_1,x_2,x_3)$, $\rho, l \in \N$,

\noindent b) $P=x_1^{\mu}x_3^kM(x_1,x_2,x_3)$, $\mu, k\in \N$,  

\noindent c) $P=x_1^{\eta}N(x_1,x_2,x_3)$, $\eta\in \N$  and $M_{\pm}\, .\, x_1=s_1 x_1$,\\
where $L,M,N\in \C[x_1,x_2,x_3]$.

\end{prop}

\begin{proof}

Since $P$ is invariant under the action of $M_+$ and \eqref{condition-1} holds, all different monomials entering in $P$ contain a positive degree of $x_1$. Therefore,  $P$ factorizes as follows
\begin{equation}
P=x_1^{\rho}\, P_1(x_1,x_2,x_3),\quad \rho\in \N\,,
\end{equation}
where $P_1$ is a  homogeneous polynomial not divisible by $x_1$.

 Let $\mathcal{E}$  be the set  of all pairwise non-colliniear forms dividing $P$. In particular we have  already $x_1\in \mathcal{E}$. We denote 
$\mathcal{E}_0\subset \C P^2$ the set of directions of all elements from  $\mathcal{E}$ equipped  with the naturally  defined  action of $G$. 
The $G$-invariance of $P$  implies  clearly  the invariance  of   directions from $\mathcal{E}_0$ under the action of $G$.

Exchanging the roles of $M_+$ and $M_-$, one finds $e\in \mathcal{E}$  -- the eigenform of $M_-$, $M_-\,.\, e=s_1e$. Let $e=e_1x_1+e_2x_2+e_3x_3$.  Then, considering the orbit $M_+^n\, .\, e$, $n\geq 1$ and using  \eqref{condition-1} together with  $\mathrm{card}\, \mathcal{E}_0<\infty$,  one sees that if $e_1\neq 0$ then $e_2=e_3=0$ i.e. $e=x_1$ (case c).  If  $e_1=0$, then either $e_2e_3= 0$ and so 
$e\in \{ x_2,x_3\}$ (cases a-b)  or $e_2e_3 \neq 0$ so that   $\exists \, m \in \N$ such that $s_2^m=s_3^m$. The last equation  implies $\lambda_2-\lambda_3\in \R$ which is impossible according to \eqref{CONDITION}. This concluded the proof.

\end{proof}

\begin{prop} \label{permut1}
The monodromy transformations $M_+$ and $M_-$, taken  in any basis, are permutationally  conjugated i.e. $\exists \, C\in GL(3,\C)$ such that 
\begin{equation}\label{permut}
C M_{+}C^{-1}=M_{-}, \quad C M_{-}C^{-1}=M_{+}\,.
\end{equation}
\end{prop}

\begin{proof}
We take a basepoint  $e\in\mathcal R$  on the positive imaginary axis $\mathrm{Im} \,p >0$. Let $\Sigma(p)$, $\Sigma(e)=\Id$ be the normalised fundamental matrix solution of \eqref{NVE1} and let $G$ be the corresponding monodromy group.  Since the equations   \eqref{NVE1}   are invariant under the change of time $p\to -p$,   $\tilde \Sigma(p)=\Sigma(-p)$ is again a fundamental matrix solutions of \eqref{NVE1}. Let $\Gamma_1$, $\Gamma_2$ be two loops  starting from  $e$ and going around the singularities   $p=\alpha$ and $p=-\alpha$, respectively.  We define the  loops $\tilde \Gamma_{i}=-\Gamma_{i}$, $i=1,2$ (symmetric to $\Gamma_{1,2}$ with respect to origin), starting from  the point $\tilde e=-e$ and having the same orientation as $\Gamma_{1,2}$. Obviously, $\tilde \Sigma (\tilde e)=\Id$ and  we can define the monodromy group $\tilde G$ using $\tilde \Sigma$ in the usual way.  One sees  that
$\tilde \Gamma_1$, $\tilde \Gamma_2$ define  now   the  monodromy transformations around $p=-\alpha$ and $p=\alpha$,  respectively. The result follows then from the fact that $G$ and $\tilde G$ are always conjugated.
\end{proof}

\begin{rem}
It follows from the previous proposition that if $C^2\neq \Id$ then it is a centraliser of $G$ in $GL(3,\C)$. 

\end{rem}

\begin{prop}\label{DIAG}
Let us assume that the monodromy group $G$ of the normal variational equations \eqref{NVE1} has a polynomial homogeneous invariant and that the conditions \eqref{condition2}, \eqref{CONDITION} hold. Then $G$ is diagonalizable. 
\end{prop}

\begin{proof}

One first assumes that   a)-b)  from Proposition \ref{pi} hold.  

Then we have three possible cases:  $\mathcal{E}=\{x_1,x_2\}$ (A), $\mathcal{E}=\{x_1,x_3\}$ (B) or  $\mathcal{E}=\{x_1,x_2,x_3\}$ (C).

 The group $G$ acts  on $\mathcal{E}_0$ by permutations.  Hence, exchanging  if necessary $x_2$ and $x_3$,  in A-B  we can put 
\begin{equation} \label{345}
M_+=\diag(s_1,s_2,s_3), \,
M_-=
\left[
\begin{array}{lll}
\sigma_1 & 0 & 0\\
0 &  \sigma_2 & 0\\
a & b & \sigma_3 
\end{array}
\right], \,  \sigma_i\in \Spectr(M_+), \, a,b\in \C\,.
\end{equation}

Indeed, $\mathcal{E}$ always contains an eigenform of $M_-$ corresponding to the stable eigenvalue $s_1\in \Spectr(M_+)=\Spectr(M_-)$. In particular $\sigma_1=s_1$ or $\sigma_2=s_1$.  

Using \eqref{345} one can calculate the unique (modulo multiplication by a constant)  nonsingular matrix $T$ such that $T^{-1} \, M_- \, T=\diag(s_1,s_2,s_3)$. Then, as follows from  Proposition  \ref{permut1}, 
$T^{-1}\, M_+\, T=M_-$. The last equation, together with \eqref{CONDITION} and some  elementary calculations  implies $a=b=0$ that  proves the result.
 
If the case C holds,  then the  matrix $M_- \in SL(3,\C)$ is either diagonal, so the proposition is proved, or is a permutation of the eigendirections of $M_+$ fixing one of them. In this last case, $\Spectr(M_-)$ will contain necessarily  a pair of eigenvalues with equal absolute values. This contradicts to \eqref{CONDITION}.

We consider now  the case c) from Proposition \ref{pi}.  One puts $P$ into the form

\begin{equation}\label{432}
P=x_1^{\rho}(x_1P_1(x_1,x_2,x_3)+P_2(x_2,x_3))\,,
\end{equation}
where $P_1$, $P_2$ are homogeneous polynomials and $P_2\neq \const$ since $|s_1|<1$.

Let $\tilde M_+=\diag(s_2,s_3)$ and let $\tilde M_-=(m_{ij})_{2\leq i,j\leq 3}$ denotes the restriction of the linear operator $M_-$  to the $x_2,x_3$-plane.

It is clear from \eqref{432}  that $P_2$ is  a polynomial semi-invariant for both $\tilde M_+$ and $\tilde M_+$.  If  $P_2$ contains two different monomials $x_2^n x_3^m$ and $x_2^p x_3^q$,  $n+m=p+q=\deg P-\rho$ then $s_2^r=s_3^r$ for $r=n-p=q-m\neq 0$ that contradicts  to \eqref{CONDITION}. 

If    $x_2^N$ (resp. $x_3^N$)  is the only monomial  entering in $P_2$, then $x_2$ (resp. $x_3$) is the eigenform of $\tilde M_-$. 

Thus, exchanging  if necessary $x_2$ and $x_3$, it is sufficient  to consider the case
\begin{equation}
M_+=\diag(s_1,s_2,s_3),\quad M_-
=
\left [ 
\begin{array}{lll}
s_1& 0 &0\\
a& m_{22}&0\\
b&m_{32}&m_{33}
\end{array}
\right]\,. 
\end{equation} 

With help of \eqref{permut}, \eqref{CONDITION} and some elementary algebraic computations one shows  that $m_{22}=s_2$, $m_{33}=s_3$.  

We introduce the matrices $U=M_+^{-1}M_-$ and $K=U-\Id$.    One verifies  that  $\Spectr(U)=\{1,1,1 \}$ and that  if $\mathrm{rang}(K)=2$  then $M_-$ does not have any polynomial invariants.  That can be done by transforming $U$ to its Jordan form.  If  $\mathrm{rang}(K)=0$ then $M_-=M_+$ and the proposition is proved.  The condition $\mathrm{rang}(I-\mathrm{Id})=1$ implies in turn:  $a=0$ (i),   $b=m_{32}=0$ (ii) or $m_{32}=0$ (iii). In the cases (i), (ii), $M_-$ has the  eigenform equal to $x_2$ or $x_3$ that corresponds to the case considered before.  If (iii) holds, it is sufficient to apply again  the conjugacy conditions \eqref{permut} to obtain a contradiction with \eqref{CONDITION}.

 Finally, we consider the case then $P_2=c\, x_2^t x_3^l$, $t,l\in N$, $c\in \C^*$.  Then $\tilde M_-$ either preserves or permutes $x_2,x_3$-eigendirections of $\tilde M_+$. In the first case  the proceed as above. In the second one,  it is sufficient to verify by a direct computation  that $\Spectr(M_-)$ will contain in this case a pair of eigenvalues  $\pm s$, $s\in \C^*$   so that the condition \eqref{CONDITION} is violated.  The proof of Proposition \ref{DIAG} is finished.

\end{proof}

The next proposition shows that in our case the Fuchsian system \eqref{NVE1} never has  a diagonal monodromy group.

\begin{prop} \label{diagonal}
Under the conditions \eqref{condition2} and \eqref{CONDITION} the monodromy group $G$ of the normal variational equations \eqref{NVE1} is not  diagonalizable.
\end{prop}

\begin{proof}

Let us assume that $G$ is diagonalizable:  $M_+=\diag(s_1,s_2,s_3)$, $M_-=\diag(s_{i_1},s_{i_2},s_{i_1})$,  $s_k=e^{2 \pi i \lambda_k}$,  $k=1,2,3$  where $(i_1,i_2,i_3)$ is a certain permutation  of $(1,2,3)$.

  That implies existence of  three independent solutions of \eqref{NVE1} of the form

 \begin{equation} \label{twodiag0}
 X_k(p)=(p-\alpha)^{\lambda_k}(p+\alpha)^{\lambda_{i_k}}N_k(p), \quad k=1,2,3\,,
 \end{equation}
 with vector functions $N_k$ rational on $p$.

Let $Y_i$, $i=1,2,3$ be three linearly independent  eigenvectors of $A$ with corresponding eigenvalues $\lambda_i$ (we remind that $\lambda_i$ are pairwise different  in view of  \eqref{CONDITION}).
One deduces from  \eqref{twodiag0} the following formulas containing  the rational vector functions $R_k$
 \begin{equation} \label{twodiag}
 X_k=(p-\alpha)^{\lambda_{g_k}}(p+\alpha)^{\lambda_{u_k}}R_k(p), \quad k=1,2,3
 \end{equation}
with  ${g_k}$ ${u_k}\in \{1,2,3 \}$  and where now $R_k=Y_{g_k}+Y^{+}_{g_k}(p-\alpha)+\cdots$  in the neighborhood of $p=\alpha$ and   $R_k=Y_{u_k}+Y^{-}_{g_k}(p+\alpha)+\cdots$ in the neighborhood of $p=-\alpha$, ($k=1,2,3$).  Since the system  \eqref{NVE1} is invariant under the change $p\to -p$,   $X_k(-p)$  is also its solution.  It yields, together with condition \eqref{CONDITION},    $g_k=u_k$   for  $k=1,2,3$.

 Under the same  condition \eqref{CONDITION}, as  a simple argument shows,   $Y_{g_1}$, $Y_{g_2}$, $Y_{g_3}$ must be  pairwise different.  Since $M_+=M_-$ and $M_+M_-M_{\infty}=\mathrm{Id}$, the similar property holds for the point $p=\infty$.

 Thus, one  can represent   $X_k$  as below
 
 \begin{equation} \label{twodiag1}
 X_k=(p-\alpha)^{\lambda_{k}}(p+\alpha)^{\lambda_{k}}p^{n_k}P_k(p), \quad k=1,2,3\,,
 \end{equation}
where $n_k\in \Spectr(B)=\{0,-1,-1\}$; $P_k$ are polynomial vector functions   such that  $P_k(\alpha)=P_k(-\alpha)=Y_k$   for $k=1,2,3$  and   $P_k(0)$  are   eigenvectors of  $B$.

Let $D_k$ be the order of $X_k$ at infinity. One finds from \eqref{twodiag1}: $D_k=-2\lambda_k-n_k-d_k$ so that $\sum D_k=-2\sum \lambda_k -\sum n_k-\sum d_k\geq -2\cdot (-1)-(0)-\sum d_k$ since $\sum n_k \leq 0$.  Otherwise, we know that  $\sum D_k=4$ (see Proposition \ref{propos1}) and hence  $\sum d_k\leq 2$.   Thus, at least one of the vectors $P_k$ is constant  and, as easy  follows from  substitution of \eqref{twodiag1} into \eqref{NVE1},   is a common eigenvector of $A$ and $B$. 
This is impossible according  to Proposition  \ref{propos2} and the proof is achieved.
\end{proof}

In view of Lemma \ref{MainProp} 
the above results  can now be summarized  
by the  main Theorem \ref{MainTh} .

\section{Conclusion}
The difficulty of the rattleback  problem is due to   its non-hamiltonian nature. That explains    the non-trivial  structure  of  the monodromy group studied in the previous sections.
The  technical problem in applying  our non-integrability Theorem \ref{MainTh}  are   the  hyperbolicity  conditions \eqref{CONDITION}. Of course, one can write these restrictions directly  using Cardano formulas that will lead to quite complicated expressions.   It may be interesting to consider a concrete example.   Let  $I_1=0.5$, $I_2=0.6$, $I_3=0.8$, $\delta=1.3$, $m=1$, $g=1$, $b_1=1$, $b_2=2$, $b_3=3$. Then the eigenvalues of the characteristic polynomial \eqref{CHARAC} are
$$
\lambda_1=-0.365-0.858 i,\quad \lambda_2= -0.435+0.963i, \quad \lambda_3=-0.200-0.106i\,,
$$
and  the conditions \eqref{CONDITION} are obviously satisfied.

We believe that our non-integrability conditions can be  strengthened. Thus, the further detailed analysis  of  the characteristic polynomial  \eqref{CHARAC} is needed. 

We note that in the  heavy rigid body case,  the characteristic polynomial \eqref{CHARAC} always has a real root  as seen from   \eqref{eh}. Indeed, it corresponds to existence of the fourth polynomial first integral  \eqref{4int} of  the Euler-Poisson equations.
In the rattleback case, the interesting remaining  problem is   the existence of new meromorphic first integrals. Our  Theorem \ref{main:lem} answers this question only  then  the variational equations \eqref{NVE1} have logarithmic singularities.   Our  intention to avoid the study of  the  Zariski closure of the monodromy group $G$  was twofold. Firstly, that makes the proofs self-contained and quite elementary.
 Secondly, we would like to underline the importance of the symmetry conditions \eqref{permut}, coming from the mechanical context  of the problem and which simplify greatly the non-integrability analysis.

\vspace{0.5cm}

\noindent { \bf  ACKNOWLEDGMENT}\\

The authors gratefully acknowledge support by the scientific program ALLIANCE  
of  the British Council and the French Foreign Affairs Ministry, Project  No. 12106TD.

\section{Appendix}

The matrix $N$ is given by 
\begin{equation}
N = \begin{pmatrix}
0 & * & * \\
\imag \frac{p^2 - \alpha}{\beta} & 0 & p \\
- \frac{ p^2 - \alpha}{\beta} & -p & 0 
\end{pmatrix} \,.
\end{equation}

The matrix $M$ is given by 
\begin{equation}
M = \begin{pmatrix}
a_1 p & a_2 p & c_1 + c_2 \frac{\alpha}{p^2-\alpha} \\
a_3 p  & -a_1 p & c_3 + c_4 \frac{\alpha}{p^2-\alpha} \\
-\imag\frac{p^2-\alpha}{\beta} & \frac{p^2-\alpha}{\beta} & 0
\end{pmatrix} \,.
\end{equation}
The coefficients $a_i$, $c_i$ are related by one quadratic equation, which is equivalent
to \eqref{coefrelation}.
They are given by 
\begin{subequations}
\begin{align}
\theta_d a_1 & = 2\frac{mg}{\beta}(\rho_1 - \rho_2) (\Sigma_{12} - \imag m \rho_1 \rho_2)
  + (\Sigma_{12} - \imag m \rho_1 \rho_2)( \Sigma_{11} + \Sigma_{22} + m(\rho_1^2 - \rho_2^2))
 \\
\theta_d a_2 & = 2\frac{mg}{\beta}(\rho_1 - \rho_2) (\Sigma_{22} + m \rho_1^2) 
\\ & \quad \nonumber
   + (\Sigma_{12} - \imag \Sigma_{22} - \imag m \rho_1( \rho_1 - \rho_2)
      (\Sigma_{12} + \imag \Sigma_{22} + \imag m \rho_1( \rho_1 - \rho_2)
\\
\theta_d a_3 & = 2\frac{mg}{\beta}(\rho_1 - \rho_2) (\Sigma_{11} - m \rho_2^2) -
\\ & \quad \nonumber
   -  (\Sigma_{11} - \imag \Sigma_{12}  - m \rho_2(\rho_1 - \rho_2) )
      (\Sigma_{11} + \imag \Sigma_{12}  + m \rho_2(\rho_1 - \rho_2) )
\\
2 \theta_d c_1 & = 
2 mg (\Sigma_{12} (\rho_1 - \rho_3) + \imag \Sigma_{22} (\rho_2 - \rho_3) +
                     \imag m \rho_1\rho_3(\rho_1 - \rho_2))  
\\                     & \quad \nonumber
- m \beta(\Sigma_{12} (2\rho_1 - \rho_2) - \imag\Sigma_{22} (\rho_1 - \rho_2)) 
+ \imag m^2 \beta\rho_1\rho_3(\rho_1^2 - \rho_2^2)\\
2 \theta_d c_3 & = 
2 mg (\Sigma_{11} (\rho_1 - \rho_3) + \imag \Sigma_{12} (\rho_2 - \rho_3) -
                     m \rho_2\rho_3(\rho_1 - \rho_2))  
\\                     & \quad \nonumber
- m \beta(\Sigma_{11} (2\rho_1 - \rho_2) - \imag\Sigma_{12} (\rho_1 - \rho_2)) 
- m^2 \beta\rho_2\rho_3(\rho_1^2 - \rho_2^2) \\
\theta_d c_2 & = -m \beta \rho_3 ( \Sigma_{12} \rho_1 + \imag \Sigma_{22} \rho_2) \\
\theta_d c_4 & = m \beta \rho_3 ( \Sigma_{11} \rho_1 + \imag \Sigma_{12} \rho_2) 
\end{align}
\end{subequations}

The matrix $B$ is given by 
\[
B = \begin{pmatrix}
-1 & 0 & -2\imag(c_1 - c_2) \\
0 & -1 & -2\imag(c_3 - c_4)\\
0 & 0 & 0 
\end{pmatrix} \,.
\]

The matrix $A$ is given by 
\[
A = \begin{pmatrix}
-\imag a_1 & -\imag a_2 & -\imag  c_2 \\
-\imag a_3 & \imag a_1 & -\imag c_4 \\
-1/\beta & -\imag/\beta & -1 
\end{pmatrix} \,.
\]

Rigid body case: 
\[
    r_1 = -\rho_1, r_2 = -\imag \rho_2, r_3 = \rho_3 = 0, \quad
    m = 0, mg \not = 0, \quad
    \Sigma_{12} = \delta = 0\,.
\]

\end{document}